\documentclass[12pt]{article}

\usepackage{amssymb,latexsym}

\newcommand{\ZZ}{\mathbb{Z}}
\newcommand{\NN}{\mathbb{N}}
\newcommand{\RR}{\mathbb{R}}

\newtheorem{thm}{Theorem}

\newtheorem{lem}{Lemma}

\newtheorem{pro}{Proposition}
\newtheorem{rem}{Remark}
\newtheorem{ques}{Questions}

\title{Always finite entropy and Lyapunov exponents of two-dimensional cellular automata}

\author{Pierre Tisseur\thanks{email: tisseur@genopole.cnrs.fr}
\\Laboratoire G\'enome et Informatique,\\
               Universit\'e d'Evry,
               Tour Evry 2,\\
               523 Place des Terrasses de l'Agora,\\
               91034 Evry Cedex  France}
\date{}

\begin{document}
\maketitle
\begin{abstract}
Given a new definition for the entropy of a cellular automata acting on a two-dimensional space, we propose an
inequality between the entropy of
the  shift on a two-dimensional lattice and some angular analog of Lyapunov exponents.
\end{abstract}

\section{Introduction}

A two-dimensional cellular automaton (CA) is a discrete mathematical idealization of a two-dimensional spacetime
physical system.  The space, consists of a discrete, infinite two-dimensional lattice with the property that each
site can take a finite number of different values. The action of the cellular automaton on this space is
the change at each time step of each values of the lattice only taking in consideration the values in a neighbourhood and
applying a local rule.

Since this automaton act on a two-dimensional lattice, this map has been used as model
in many areas (Physics, Biology ...) (see Wolfram\cite{Wo86}) and seems to be more usefull in pratice than the one-dimensional case.
Nevertheless there is very few mathematical results in the two-dimensional case except the work of
Willson \cite{Will}, Margara \cite{FMM98} about ergodicity and density of periodic points of linear
two-dimensional (CA).

The first cellular automaton defined by Von Neumann and Ulam was two-dimensional one for theoritical
self-reproducting biological systems as the well knwon game of life defined by J. Conway.

In \cite{Sh91} Shereshevsky  defined the first
Lyapunov exponents for one-dimen\-sional cellular automata and establish and inequality with the
standart metric entropy.
Then he asked for a generalisation of these results in higher dimensions and qualify this extension as a
challenging problem.
This is the topic of this paper : first  defined a natural equivalent of metric entropy which is finite
for all two-dimensional (CA) and try to rely this value with a generalisation of (CA) Lyapunov
exponents in the two-dimensional case.

Entropy is an isomorphism invariant for dynamical systems see (Denker \cite{DGS75})and can be seen as a measure of the disorder of the
dynamical system. It is always finite for a one-dimensional (CA).
In the two-dimensional case the metric entropy is in general not finite  when the shift metric entropy
is positive. M. Shereshevsky postulate that a two-dimensional (CA) have an entropy  equal to zero or infinite.
Some results shows that for additive one the conjecture is true.

So we propose here a new  measure of complexity for two-dimensional (CA) : the Always Finite Entropy (AFE) of
a two-dimensional cellular automaton $F$ denote by $h_\mu^{2D}(F)$.
In section 3 we show by examples that the (AFE) give equivalent values of entropies of
one-dimensional cellular automaton rules embeded in a two-dimensional lattice.

Then like in the one-dimensional case (see Shereshevsky \cite{Sh91}, Tisseur \cite{ti2000} ) we defined analog of Lyapunov
exponent of differential dynamical systems for (CA).
These Lyapunov exponents can be seen as the speed of propagation of perturbation in a particular direction.
These Lyapunov exponent can be seen more like a generalisation of those defined in the one-dimensional case
than a direct transciption of the differential case.

Hence using some generalization of Shannon theorem we obtain an inequality between the always finite entropy
of the (CA), the two-dimensional entropy $h_\mu (\sigma_1,\sigma_2)$ and a map which depends
on directional Lyapunov exponents :
$$
h_\mu^{2D}(F)\le h_\mu (\sigma_1,\sigma_2)\times \left(\int \frac{\lambda^2(\theta)}{2}d\theta+
\sqrt(2)\int \lambda(\theta )d\theta +
\frac{\pi}{2}\right)
$$
It seems that this  result can be generalised in higher dimensions.

\section{Definitions}

\subsection{Symbolics systems and  cellular automata}

Let $A$ be a finite set or alphabet. Denote by $A^{**}$ the set of all
two-dimensional concatenations of letters in $A$.
The finite concatenations are called blocks.
The number of letters of a block $u\in A^{**}$ is denoted by  $\vert u\vert$.
The set of infinite lattice
$x=(x_{(i,j)})_{(i,j)\in\ZZ^2}$ is denoted by $A^{\ZZ^2}$. A point $x\in A^{\ZZ^2}$ is called a
configuration. For each integers $s$ and $t$ and each block $u$, we
call cylinder the set $[u]_{(s,t)}=\{y\in A^{\ZZ^2} : y_{(s+i,t+j)}=u_{(i,j)}\}$.
For $i_1\le i_2$ and $j_1\le j_2$ in $\ZZ$ we denote by $x((i_1,i_2),(j_1,j_2))$ the rectangular block
 $x_{(i_1,j_1)}\ldots x_{(i_1,j_2)};\; x_{(i_1+1,j_1)}\ldots x_{(i_1+1,j_2)};\ldots\ldots
 x_{(i_2,j_1)}\ldots x_{(i_2,j_2)}$.
We endow $A^{\ZZ^2}$ with the product topology.
For this topology $A^{\ZZ^2}$ is a compact metric space.
A metric compatible with this topology can be defined by the distance
$d(x,y)=2^{-t}$ where $t=\min\{s=\sqrt{i^2+j^2}\; \vert \,\mbox{ such that }
x_{(i,j)}\ne y_{(i,j)}\}$.
The horizontal shift $\sigma_1 \colon A^{\ZZ^2}\to A^{\ZZ^2}$ is defined by :
$\sigma_1 (x)=(x_{(i+1,j})_{(i,j)\in \ZZ^2}$.
The vertical shift $\sigma_2 \colon A^{\ZZ^2}\to A^{\ZZ^2}$ is defined by :
$\sigma_1 (x)=(x_{(i,j+1})_{(i,j)\in \ZZ^2}$.

The dynamical system $(A^{\ZZ^2} ,\sigma_1,\sigma_2 )$ is called the full shift. A subshift $X$
is a closed shift-invariant subset $X$ of $A^{\ZZ^2}$ endowed with the shifts
$\sigma_1$ and $\sigma_2$.
It is possible to identify $(X,\sigma_1,\sigma_2 )$ with the set $X$.

 If $\alpha =\{A_1,\ldots ,\,
A_n\}$ and $\beta =\{B_1,\ldots ,\, B_m\}$ are two partitions of $A^{\ZZ^2}$ denote
by $\alpha \vee \beta$ the partition $\{A_i\cap A_j\, i=1\ldots n, \,\,
j=1,\ldots ,\, m\}$.
The metric entropy $h_\mu (T)$ of a transformation $T$
is an isomorphism invariant between two $\mu$-preserving
transformations. Put $H_\mu (\alpha ) = \sum_{A\in\alpha}\mu (A)\log \mu (A)$ where $\alpha$ is a finite
partition of the space.
The entropy of the partition $\alpha$ is defined as $h_\mu (\alpha ) =
\lim_{n\to\infty}1/nH_\mu (\vee_{i=0}^{n-1}T^{-i}\alpha )$ and the entropy of $(X,T,\mu )$ as
$\sup_\alpha h_\mu (\alpha )$.

 A two-dimensional cellular automaton  is a continuous self-map $F$ on
$A^{\ZZ^2}$ commuting with the two-sided shift. We can extend the Curtis-Hedlund-Lyndon theorem \cite{DGS75}
an state  that for every two-dimensional cellular automaton $F$
there exist an integer $r$ and a
block map $f$ from $A^{r^2}$ to $A$ such that:

$F(x)_{(i,j)}=f(x_{i-r,j-r},\ldots ,x_{(i,j)},\ldots ,x_{(i+r,j+r)}).$

As in the one-dimensional case we call  radius of the cellular
automaton  $F$ the integer $r$ which appears in the definition of the associated local rule.

\subsection{Entropy of two-sided shift}

Recall that $\sigma_1$ and $\sigma_2$ are respectively the horizontal and the vertical shift acting
on the space $A^{\ZZ^2}$.
Let $\alpha_1$ be the partition of $A^{\ZZ^2}$ by the central coordinate $(0,0)$ and
$(A_n)_{n\in\NN}$ be a sequence of finite sets of $\ZZ^2$.
 We denote by $\vert A_n\vert$ the number of element of $A_n$.
Define $\alpha_n^{\sigma_1,\sigma_2}= \vee_{(i,j)\in A_n} \sigma_1^{-i}\sigma_2^{-j}\alpha_1$ and
$\alpha_n^{\sigma_1,\sigma_2}(x)$ as the element of the partition $\alpha_n^{\sigma_1,\sigma_2}$ which
contains the point $x\in A^{\ZZ^2}$.

In \cite{Or83} Ornstein and Weiss give an extension of the Shannon McMillan theorem for
 a class of amenable groups. These results can be used for $\ZZ^2$ actions if the sequence of
 the finite partitions verifies some special conditions.
First the set $A_n$ must be averaging set, this mean that
$$
\lim_{n\to\infty}\vert \sigma_i (A_n)\Delta A_n\vert /\vert A_n\vert =0, \;\;\;
\mbox{ with i=1 or i=2}
$$
and $\Delta$ denotes the symmetric difference of two sets.

Then to be special averaging,  the sequence $(A_n)_{n\in\NN}$ have to satisfy
$A_1\subset A_2\subset \ldots A_n$.

Finaly if  $A_n$ is a special averaging sequence,
  for almost all points $x\in A^{\ZZ^2}$ (see \cite{Or83}) one has

$$
h_\mu (\sigma_1,\sigma_2,\alpha_1 )=\lim_{n\to\infty}\frac{1}{\vert A_n\vert}
-\log \mu \left( \alpha_n^{\sigma_1,\sigma_2}(x)  \right).
$$
Remark that $h_\mu (\sigma_1,\sigma_2,\alpha_1 )=h_\mu (\sigma_1,\sigma_2 )$ because $\alpha_1$ is
a generating partition for $\sigma_1,\sigma_2$.

We are going to use this result in the case of sequences $(A_n)_{n\in\NN}$ which are not
rectangles (in this case we have the equality by definition ) and show that
$$
\int_{A^{\ZZ^2}}\lim_{n\to\infty}\frac{-\log \mu \left(\alpha_n^{\sigma_1,\sigma_2}(x)\right)}
{\vert A_n\vert}d\mu (x)=h_\mu (\sigma_1,\sigma_2).
$$

\subsection{Always finite entropy of a two-dimensional CA}

Let $F$ be a two-dimensional cellular automata and $\mu$ a shift ergodic and $F$ invariant
measure. Let $\alpha$ be a finite partition of $A^{\ZZ^2}$.
Let $h_\mu (F,\alpha )$ be the entropy of $F$ with respect of the
finite partition $\alpha$.
If $h_\mu(\sigma_1,\sigma_2)>0$
then in general it  exists a sequence of partition $\alpha_n$ such that $h_\mu (F,\alpha_n )$ goes to
infinty (see the examples below).

For this reason we propose a new kind of entropy for a two-dimensional cellular automata $F$, the always
finite entropy (AFE) denoted by $h^{2D}_\mu (F)$.
For this definition we need to defined a sequence of finite partitions $\alpha^F_n$.

Let $B$ be a two-dimensional squarred block that is to say a double finite squarred sequence of letters
 $B=(B_{(i,i)})_{(0\le i\le n; 0\le j\le n)}$.
We denote by $[B]^{(n,n)}_{(0,0)}$ the cylinder which is the set of all the points
 $y$ such that  $y_{(i,j)}=B_{(i,j)}$ with $0\le i\le n$ and $0\le j\le n$.
 We call $n$-squarred cylinders blocks the cylinders $[B]^{(n,n)}_{(0,0)}$.

For each positive integer denote by $\alpha_n$ the partition of $A^{\ZZ^2}$ into $\vert A\vert ^{n^2}$
 $n$-squarred cylinders blocks $[B]^{(n,n)}_{(0,0)}$.

Then define
$$
h^{2D}_\mu (F)=\liminf_{n\to\infty}\frac{1}{n}h_\mu (F,\alpha_n).
$$

Later we will see (main theorem ) why this entropy is always finite.

\begin{ques}
If we only ask to fix a number of $k\times n^2$ coordinates ($k\in\NN^*$), in a partition $\alpha_n^*$, is
$\liminf_{p\to\infty}h_\mu (F,\alpha_n^*)=h_\mu^{2D}(F)$?
\end{ques}

\begin{ques}
If what case  $\left(\frac{1}{n}h_\mu (F,\alpha_n)\right)_{n\in\NN}$ is
a converging sequence?
\end{ques}


\subsection{Directional Lyapunov exponents}

An almost line $L^\theta_{(i,j)}$ is a doubly infinite sequence of coordinates
$(u_n,v_n)_{n\in\ZZ}$ with $(u_0,v_0)=(i,j)$. Let $E(x)$ be the integer part of a real $x$.
Consider the continuous space $\RR^2$ and a line which pass by the coordinate $(i,j)$ and make an angle
$\theta$ with the vertical line. This line is the set of points $(x,y)$. Starting from $(x,y)=(i,j)$
and making grow the variable $x$ the succession of value $(E(x),E(y))$ gives the sequence $(u_n,v_n)_{n\in\NN}$
of the positive part of the almost line $L^\theta_{(i,j)}$.
  We obtain the negative coordinate by decreasing the variable $x$ and taking again the
 sequence $(E(x),E(y))$.

Let $W_{(i,j)}^\theta\subset\ZZ^2$ be the half plane of the lattice consisting of all the coordinates
situated in the halph plane delimited by the almost line $L^\theta _{(i,j)}$ and which contains the
central coordinate $(0,0)$.
For all point $x\in A^{\ZZ^2}$ note $W_{(i,j)}^\theta (x)$ the set of all points $y$ such that
$y_{(a,b)}=x_{(a,b)}$ for $(a,b)\in W_{(i,j)}^\theta$.

Let define the propagation map of information  : $\Lambda_n^F(\theta )$ during $n$ iterations
in the direction $\theta$ for a
two-dimensional cellular automaton $F$ :
$$
\Lambda_n^F(\theta)(x)=\min \{s=\sqrt{i^2+j^2}\vert
F^k \left(W_{(i,j)}^\theta (x)\right) \subset W_{(0,0)}^\theta (F^k(x)) \vert 1\le k\le n\}.
$$
Put $\Lambda_n^F(\theta)=\max\{\Lambda_n^F(\theta)(x)\vert x\in A^{\ZZ^2}\}$.
Remark that for each real $\theta$  there exist also a couple
$(i,j)=G_nF(\theta)$ of positive integers such
that $\sqrt{i^2+j^2}=\Lambda^F_n(\theta )$.

\begin{lem}
For all real $\theta\in [0,2\pi]$ the limit
$\lim_{n\to\infty}\frac{\Lambda^F_n(\theta )}{n}$ exits.
We denote by $\lambda (\theta )$ these limits.
\end{lem}
{\it Proof}

First prove that for each couple of integers $m$, $n$ and point $x\in A^{\ZZ^2}$, one has
$\Lambda_{m+n}^F (\theta )\le \Lambda_n^F(\theta )+\Lambda_m^F(\theta )$.

To simplify put $(s,t)=G_nF(\theta)$ and $(u,v)=G_nF(\theta')$,
$s\approx\lfloor\Lambda^F_n(\theta)\cos(\theta)\rfloor$ and
$t\approx\lfloor\Lambda^F_n(\theta)\sin(\theta)\rfloor$
.

From the definition of $s,t$ one has for all point $x\in A^{\ZZ^2}$
$F^n(W_{(s,t)}^\theta (x))\subset W_{(0,0)}^\theta (F^n(x))$.
Then
$$
F^{m+n}(W_{(s,t)}^\theta (x))\subset F^m(W_{(0,0)}^\theta (F^n(x)))=F^m(W_{(u,v)}^\theta
(\sigma_1^u\circ\sigma_2^v\circ F^n(x))).
$$
And by definition of $u$ and $v$ we obtain
$$
F^m(W_{(u,v)}^\theta (\sigma_1^u\circ\sigma_2^v\circ F^n(x)))\subset W_{(0,0)}^\theta (\sigma_1^u\circ
\sigma_2^v\circ F^{n+m}(x)).
$$
Then $F^{m+n}(W_{(s+u,t+v)}^\theta(\sigma_1^u\circ\sigma_2^v(x)))\subset W_{(0,0)}^\theta (\sigma_1^u\circ
\sigma_2^v\circ F^{m+n}(x))$,
this implies that $F^{m+n}(W_{(s+u,t+v)}^\theta(x))\subset W_{(0,0)}^\theta (F^{m+n}(x)).$

Finaly we can conclude arguing that $\Lambda^F_{n+m}(\theta)\le \sqrt{(s+u)^2+(t+v)^2}=
\sqrt{s^2+t^2}+\sqrt{u^2+v^2}=\Lambda^F_n(\theta )+\Lambda^F_m(\theta )$.

If we fix $\theta$, the numeric sequence $\frac{\Lambda^F_n(\theta )}{n}$ is a subadditive sequence
an so converge to $\inf_n \frac{\Lambda^F_n(\theta )}{n}$.

\hfill$\Box$

\subsection{Lyapunov exponent surface}
The computation of entropy of $\ZZ^2$ action using the Shannon-Breiman-Mc\-Millan Theorem need
to know the number of coordinates fixed for an element of the partition $\alpha_n^{n,F}$ (see section
2.2).
This number of elements can be seen as the surface of a polygone on the squarre lattice.

First define the sequence of maps $G_n$ from $[0,2\pi]$ to $\ZZ^2$ such that $G_n(\theta )=
(\lceil \sqrt{2}n\cos{\theta}\rceil$ ; $\lceil \sqrt{2}n\sin{\theta}\rceil )$ if $\theta\in [0,\pi]$
and $G_n(\theta )=(0,0)$ if $\theta\in ]\pi,2\pi]$.
Let $T_n$ be the set of coordinates  of the squared lattice in the interior
of the set of coordinates $(i,j)=G_nF(\frac{2\pi k}{n})+G_n(\frac{2\pi k}{n} )$, $0\le k\le n$.
Let $T_n=\{(a,b)\in\NN^2\vert \exists (i,j)=G(\theta ) (\theta\in[0,2\pi])\; \vert
\sqrt{a^2+b^2}\le\sqrt{i^2+j^2}+G(\theta )\}$.

Remark that we decide arbitrairely to rely the steps in the variations of the angle $\theta$
($\triangle\theta =\frac{2\pi}{n}$) to define the set $T_n$ with the number of iterations
of the cellular automaton. Some small changes in the definition of $T_n$ will not change
the main result.

Let define a new surface $T_n^*$.
The surface $T_n^*$ is the intersection of all the  almost half planes
$W_{G(\frac{2\pi k}{n})+GF(\frac{2\pi k}{n})}^{\frac{2\pi k}{n}}$
 (defined in section 2.4) with $0\le k\le n$.
We call $T_n$ the Lyapunov exponent surface and $T_n^*$ the surface of common behaviour.
As it is difficult to give an expression of the surface $T_n^*$ using the Lyapunov exponents $\lambda(\theta)$
in order to establish an inequality with $h_\mu^{2D}(F)$,
we are going to show that the part of $T_n^*$ wich do not belong to $T_n$ became very small in comparaison
of $T_n$ when $n$ increase.

Let $T_n^{**}=T_n^*-(T_n^*-T_n\cap T_n^*)$.
\begin{lem}
For each point $x\in A^\ZZ$ one has
$$
\lim_{n\to\infty}\frac{\vert T_n^{**}\vert}{\vert T_n^*\vert}=1
$$

\end{lem}
{\it Proof}

Let $\vert T^{**}_n\vert$ and $\vert T^*_n\vert$ be the respective surface of the sets $T_n^{**}$ and $T_n^*$.
Let $DT_n= \vert T^*_n\vert -\vert T^{**}_n\vert$, one has $\frac{\vert T^{**}_n\vert}{\vert T^*_n\vert}=
\frac{\vert T^*_n\vert -DT_n}{\vert T^*_n\vert}$.
As the sets $T^*_n$ and $T^{**}_n$ are polygones the difference of the two surfaces is at least a
sum of $n$ surfaces of triangles. We decide to bound the surface of these triangles by surfaces of
rectangles.

For each $0\le i\le n$ we consider the rectangle $r_i$ defined by the points $p_1$
$GF(\frac{\pi}{2}-\frac{2\pi i}{n})+
G(\frac{\pi}{2}-\frac{2\pi i}{n})$; $GF(\frac{\frac{\pi}{2}2-\pi (i+1)}{n})+
G(\frac{\pi}{2}-\frac{2\pi (i+1)}{n})$, the point $p_2$ which is the intersection of the almost line
$L_{(0,0)}^{\frac{\pi}{2}-\frac{2\pi}{n}}$ and $L_{GF(\frac{2\pi (i+1)}{n})+
G(\frac{2\pi (i+1)}{n})}^{\frac{2\pi}{n}}$ and the fourth point which closed the rectangle  .

Denote $\vert r_i\vert$ the surface of this rectangle.
 If the point $p_2$ is at the left of side of $p_1$ then put $\vert r_i\vert=0$.
One has $DT_n\le \sum_{i=+}^n \vert r_i\vert$.
One  has $r_i\le l_i\times h_i$ where $l_i$ is the width and $h_i$ the lenght.
One has $l_i\le h_i\times \tan (\frac{2\pi}{n})$. 

Hence

$r_i=l_i\times h_i\le (h_i)^2\tan (\frac{2\pi}{n})$ and
$h_i=\left(\Lambda_n^F (\frac{2\pi}{n})+\sqrt{2}n\cos{2\pi i/n}\right)\tan(\frac{2\pi}{n})$ so
$r_i\le (\Lambda_n^F(\frac{2\pi}{n})+\sqrt{2}n\cos{2\pi i/n})^2(\tan(\frac{2\pi}{n}))^3$ and when $n$ goes
big enought
$r_i\le (\Lambda_n^F(\frac{2\pi}{n})+\sqrt{2}n\cos{2\pi i/n})^2(\frac{2\pi}{n})^3\approx \frac{Kn^2}{n^3}
=\frac{K}{n}$ where $K$ is a constant.
Therefore $DT_n\le n\times r_i\le \frac{K}{n}\times n=K$.
So $\frac{\vert T^{**}_n\vert}{\vert T^{*}_n\vert}\le \frac{\vert T^*_n\vert -K}{\vert T^*_n\vert}$
then as we suppose $\vert T^*_n\vert$ goes to infinity then
$\lim_{n\to\infty}\frac{\vert T_n^{**}\vert}{\vert T_n^*\vert}=1$.

\hfill$\Box$

Let $\alpha_1$ the partition by the  central coordinates defined in the section 2.3 and $\alpha_1(x)$
the element of the partition which contains de point $x$.
Denote by $\alpha_n^{\sigma ,T_n^*}(x)$ the element of the partition
$\vee_{(i,j)\in T_n^*} \sigma_1^{i}\sigma_2^{j}(\alpha_1)$ which contains de point $x$.
\begin{pro}\label{pro1}
For all $0\le k\le n$  and all $x\in A^{\ZZ^2}$ one has
$F^k\left(\alpha_n^{\sigma ,T_n^*}(x)\right)\subset \alpha_n(F^k(x))$ and
$\alpha_n^F(x)\supset\alpha_n^{\alpha ,T_n^*}(x)$.
\end{pro}
{\it Proof}

Suppose that there exists $y\in \alpha_n^{T_n^*}(x)$ such that there is some $0\le k\le n$ with
$$
F^k(y)\left((0,0);(n,n)\right)\neq F^k(x)\left((0,0);(n,n)\right) .
$$
We are going to define a sequence $(Z_l(x))_{l\in\NN}$ of subsets of $\alpha_n^{T_n^*}(x)$.
Let
$$
Z_l(x)=\{z\in X\vert z\in \alpha_n^{T_n^*}(x)\cap_{m=0}^l
W^{\frac{2\pi m}{n}}_{G_n(\frac{2\pi m}{n})+G_nF(\frac{2\pi m}{n})}(x)\}.
$$
Remark that $Z_i(x)\supset Z_{i+1}(x)$ and $Z_n(x)=\{x\}$.

From the supposition there exist $0\le l<n$ such that there exist $u\in Z_l(x)$ and $v\in Z_{l+1}(x)$
such that for all $0\le k\le n$ one has
$F^k(v)\left((0,0);(n,n)\right)=F^k(x)\left((0,0);(n,n)\right)$ and there exists $0\le k\le n$
such that $F^k(u)\left((0,0);(n,n)\right)$ $\neq F^k(x)\left((0,0);(n,n)\right)$.
As $u\in Z_l(x)$ and $v\in Z_{l+1}(x)$ then
$$
v\in W^{\frac{2\pi l+1}{n}}_{G_n(\frac{2\pi l+1}{n})+G_nF(\frac{2\pi l+1}{n})}(u).
$$
It follows that
$$
F^i\left(W^{\frac{2\pi l+1}{n}}_{G_n(\frac{2\pi l+1}{n})+G_nF(\frac{2\pi l+1}{n})}(u) \right)\backslash
\!\!\!\!\!\!\subset
W^{\frac{2\pi l+1}{n}}_{G_n(\frac{2\pi l+1}{n})}(F^i(u))
$$
and
$$
F^i\left(W^{\frac{2\pi l+1}{n}}_{G_nF(\frac{2\pi l+1}{n})}(\sigma^{-G_n(\frac{2\pi l+1}{n})} (u) \right)
\backslash\!\!\!\!\!\!\subset
W^{\frac{2\pi l+1}{n}}_{(0,0)}(F^i(\sigma^{-G_n(\frac{2\pi l+1}{n})}(u))
$$
then $\Lambda_n^F(\frac{2\pi l+1}{n})>\sqrt{i^2+j^2}$ where $(i,j)=G_nF(\frac{2\pi l+1}{n})$, which
which contradict the hypothesis, so we can conclude.

\hfill$\Box$

\section{Examples of $h_\mu^{2D}$ computation and a first inequality}

In this section we give different justification of the choice of the definition of
the Always Finite Entropy. From this section we denote by $X$ the full shift $A^{\ZZ^2}$.

\subsection{Examples}
Using tree examples we are going to show that the (AFE) is a natural extension of the standart entropy
apply in the one-dimensional case.

In the three examples we consider the two-dimensional lattices on an alphabet which contains
two letters. Let $A=\{0,1\}$ and consider the space $X=A^{\ZZ^2}$.
We endows the space $X$ with the uniform measure $\mu$. This measure is the unique measure
which gives the same weight to all cylinders with the same number of fixed coordinates.

$\bullet$ The first example ($F_1$) is  the horizontal shift named $\sigma_1$ in section 2.2 .
This (CA) is defined by the rule $[F_1(x)]_{(i,j)}=x_{(i+1,j)}$.
Using the definition of metric entropy  one has
$$
h_\mu (F,\alpha_1)=\int_X -\lim_{n\to\infty}\frac{1}{n}\log \mu (\alpha_1^{n,F_1}(x))d\mu(x).
$$
As $\mu$ is the uniform measure and $F_1$ act on each point $x$ in the same way it follows that
$h_\mu (F_1,\alpha_1)= -\lim_{n\to\infty}\frac{1}{n}\log \mu (\alpha_1^{n,F_1}(x))$.

Clearly $\alpha_1^{n,F_1}(x)=x((0,0);(n,0))$ so we have $\mu (\alpha_1^{n,F_1}(x))=2^{-n-1}$, then
 $h_\mu (F_1,\alpha_1)=\lim_{n\to\infty}\frac{1}{n}(n-1)\log (2)=\log (2)$.

Using the same methods we obtains $h_\mu (F_1,\alpha_p)=p\log (2)$ and therefore
$h_\mu^{2D}(F_1)=\lim_{n\to\infty}\frac{1}{p}p\log (2)=\log (2)$.

We remark that the value of $h_\mu^{2D}(F_1)$ is the same that the value of the standart entropy
of the shift in the one-dimensional case and is equal to the entropy of the two shift action
$h_\mu (\sigma_1,\sigma_2 )$ on $X$.

\bigskip

$\bullet$ The second example named $F_2$ is a (CA) based on  an additive rule on the right
and left first neigbhours.
This (CA) is defined by  $[F(x)]_{(i,j)}=x_{(i-1,j)}+x_{(i+1,j)}$ modulo 2.
In order to take in account the right and left in an independantly way we need to consider  partitions
 $\alpha_n$ with $n\ge 2$.
Clearly one has $\alpha_2^{n,F_2}(x)=x((-n,2);(n+1,0))$ so
$\mu (\alpha_2^{n,F_2}(x))=2^{-2(2n+2)}$ and in general $\mu (\alpha_p^{n,F_2}(x))=2^{-p(2n+2)}$
then $h_\mu^{2D}(F_2)=\lim_{p\to\infty}\frac{1}{p}2p\log (2)=2\log (2)$.

Remark that the value of $h_\mu^{2D}(F_2)$ is the same that the entropy of the simple additive (CA)
in the one-dimensional case.  

\bigskip

$\bullet$ The third example ($F_3$) is an additive rule with the upper and right first neigbhours.
This (CA) is defined by the rule $[F(x)]_{(i,j)}=x_{(i,j+1)}+x_{(i+1,j)}$ modulo 2.


Recall that the set  $\alpha_p^{n,F_3}(x)$ is the element of the partition $\vee_{i=0}^n F_3^{-i}\alpha_p$
which contains the point $x$.
First remark that for any $x\in X$ one has $\mu (\alpha_2^{0,F_3}(x))=2^{-4}$ because we have to fix 4 coordinates, the squarre
delimited by the position $(0,0)$ to $(1,1)$.

Fix $x\in X$ and choose $y$
such that $y_{\left((0,0);(1,1)\right)}$ = $x_{\left((0,0);(1,1)\right)}$;
$F_3(y)_{\left((0,0);(1,1)\right)}$ = $F_3(x)_{\left((0,0);(1,1)\right)}$.
In this case we  have to fix $y_{(0,2)}$ and  $y_{(2,0)}$ but
we have the possibility to choose
$y_{(1,2)}$ or $y_{(2,1)}$. We can not choose these two coordinates in an independantly way.
 We deduce that
$\mu (\alpha_2^{1,F_3}(x))$=$2\times 2^{-4}\times 2^{-4}=2^{-7}$.

More generaly for $\alpha_p(x)$ we have to fix $p^2$ coordinates, then in $\alpha_p^{1,F}(x)$
we must fix $p^2+2p$ coordinates $(y_{(0,p)}\ldots y_{(p-1,p)}$ and
$(y_{(p,p-1)}\ldots y_{(p,0)}$)  but there is two ways of choosing the coordinates in the angle
$y_{(p-1,p)}$ and $y_{(p,p-1)}$ so $\mu (\alpha_p^{1,F_3}(x))$=$2\times 2^{-(p^2+2p)}=2^{2(p-1)}$.

In order that $F^2_3(y)\in\alpha_p(F^2(x))$ we have to fix $2p+1$ more coordinates in $y$ with
again two possibilities which multiply the first two possibilities,  so
$\mu (\{y\vert F^2_3(y)\in\alpha_p(F^2(x))\})=2^2\times 2^{-(p^2+2p+1)}=2^{-(2p-1)}$
and
$$
\mu (\{y\vert F^2_3(y)\in\alpha_p(F^2(x));\; F_3(y);\; F(y)\in\alpha_p(F(x));\; y\in\alpha_p(x)\})
$$
$$
=2^{-(p^2)}\times 2^{-(2p-1)}\times 2^{-(2p-1)}=2^{p^2+2(p-1)}.
$$

In general,  considering all the configuration such that $F_3^i(y)\in$ $\alpha_p(F_3^i(x))$ ($0\le i\le n$),
we have to fix
 $2p+i$ more coordinates   (with $2^{i+1}$ possible choices) in order that
 $F_3^{i+1}(y)\in \alpha_p\left(F_3^{i+1}(x)\right)$.
The set $\alpha_p^{n,F_3}(x)$ is the union of $\pi_{i=1}^n 2^i$ cylinders block with a number of
$2p+i$ fixed coordinates.
  Hence in order to have $F_3^i(y)\in \alpha_p(F_3^i(x))$ with
$0\le i\le n$ we have to fix $p^2+\sum_{i=0}^{n-1} (2p+i)$ coordinates but there is
$\pi_{i=1}^n 2^i$ possible choices of coordinates.
Finaly as all the configurations have the same weight
for each $x\in X$ and one has $\mu \left(\alpha_p^{n,F_3}(x)\right)=
2^{\left(-\left(p^2+\sum_{i=0}^{n-1} (2p+i)\right)+\sum_{i=1}^n 2^i\right)}
=2^{-(p^2+(2p-1)n)}$.

 As in general $\mu (\alpha_p^{n,F_3}(x))=2^{-((2p-1)n+p^2)}$ then
$$
h_\mu (F_3,\alpha_p)=\lim_{n\to\infty}\frac{(2p-1)n+p^2}{n}\log (2)=(2p-1)\log (2)
$$
and
$h_\mu^{2D} (F)=\lim_{p\to\infty}\frac{1}{p}(2p-1)\log (2)=2\log (2)$.

Remark that $h_\mu^{2D}(F_3)=h_\mu^{2D}(F_2)$ which seems natural and coherent.

\subsection{Equivalent and upper bound maps of $h_\mu^{2D}$}
Here we define two maps which can be related to the always finite entropy. The first one is an upper bound
required to establish the
main inequality with  the Directional Lyapunov Exponent, the second one is an equivalent
definition of the (AFE).

\begin{pro}\label{firstinequality}
For each $F$ invariant measure $\mu$  one has
$$
h_\mu^{2D}(F)\le \int_X \liminf_{n\to\infty}\frac{-\log \mu (\alpha_n^{n,F}(x))}{n^2} d\mu (x).
$$
\end{pro}
{\it Proof}

Let $p_i$ be a sequence such that
$\lim_{n\to\infty}\frac{1}{p_i}h_\mu (\alpha_{p_i},F)=h_\mu^{2D}(F)$.
Given any real $\epsilon >0$, there exists an integer $I$ such that if $i\ge I$ then
$\vert h_\mu^{2D}(F)-\frac{1}{p_i}h_\mu (\alpha_{p_i},F)\vert \le \epsilon$.

From the definition of the metric entropy and the dominated convergence theorem , one has
\small
$$
h_\mu (\alpha_p,F)=\int_X \lim_{n\to\infty}\frac{-\log (\mu (\alpha_p^{n,F}(x)))}{n} d\mu (x)=
\lim_{n\to\infty}\int_X\frac{-\log (\mu (\alpha_p^{n,F}(x)))}{n} d\mu (x).
$$
\normalsize

As the sequence
$$
\left(\int_X -\log (\mu (\alpha_p^{n,F}(x)))d\mu (x)\right)_{n\in\NN}
=  \left(H(\alpha_p^{n,F})\right)_{n\in\NN}
$$
 is a subadditive sequence (cf Denker \cite{DGS75}), then
$\left(\frac{H(\alpha_p^{n,F}}{n})\right)_{n\in\NN}$ converge to
$\inf_n\frac{H(\alpha_p^{n,F})}{n}=h_\mu (\alpha_p,F)$.
Then for all $n$ and $p$ in $\NN$ one
has $\frac{1}{p}h_\mu (\alpha_p,F)\le \int_X\frac{-\log (\mu (\alpha_p^{p,F}(x)))}{p^2}d\mu (x)$
and taking $\epsilon\to 0$ we can conclude.

\hfill$\Box$

The second proposition give a new seing of the always
finite entropy.

\begin{pro}
For any $F$ invariant measure $\mu$  one has
$$
h_\mu^{2D} (F)=
\lim_{\epsilon\to 0}\liminf_{n\to\infty}\int_X \frac{-\log \left(\mu \left(\alpha_
{\lceil\epsilon\times n\rceil}^{n,F}(x)\right)
\right)}{\lceil\epsilon\times n^2\rceil}d\mu (x).
$$
\end{pro}
{\it Proof}

Let $(p_i)_{i\in\NN}$ be a subsequence such that
$\lim\frac{1}{p_i}h_\mu (F,\alpha_{p_i}^{p_i,F})=h_\mu^{2D}(F)$.

Given any real $\eta >0$, there exists an integer $I$ such that if $i\ge I$ then
$\vert h_\mu^{2D}(F)-\frac{1}{p_i}h_\mu (F,\alpha_{p_i}^{p_i,F})\vert \le \eta/3$.

Then there exist an integer $N_1$ such that for all $n\ge N_1$ one has
$$
\left\vert \frac{1}{p}h_\mu (\alpha_p ,F)-\frac{1}{p}\int_X
\frac{-\log\left( \alpha_p^n(x)\right)}{n}d\mu (x)\right\vert \le \eta /3 .
$$
Take $\epsilon <\frac{P}{N_1}$, and define $n_i$ a subsequence such that $\lceil\epsilon n_i\rceil=p_i$.
Let $L(\epsilon )=\liminf_{n_i\to\infty}\int_X \frac{-\log\left( \alpha_p^{n_i}(x)\right)}{n_i}d\mu (x)$.
Denote by $n_j$ a subsequence of $n_i$ such that $\int_X \frac{-\log\left( \alpha_p^{n_j,F}(x)\right)}{n_j}
d\mu (x)$
converge to $L(\epsilon )$.
There exist an integer $N_2$ such that if $n_j\ge N_2$ one has
$$
\left\vert \int_X \frac{-\log\left( \alpha_{\lceil\epsilon n_j\rceil}^{n_j,F}(x)\right)}{\lceil{n_j}^2\epsilon\rceil}d\mu (x)
-L(\epsilon )\right\vert\le \eta /3.
$$
If $N_1\ge N_2$ take $N =N_1$ and if $N_1\le N_2$ put $P=\lceil \epsilon N_2\rceil$ in the first condition
 and $N=N_2$.
As the sequence $\int_X
\frac{-\log\left( \alpha_p^n(x)\right)}{n}d\mu (x)$ is 	a decreasing sequence we  obtain that for all
$\epsilon'\le \epsilon $ and all $n_j\ge N$
$$
\left\vert h_\mu^{2D}(F)-\int_X \frac{-\log\left( \alpha_{\lceil\epsilon n_j\rceil}^{n_j,F}(x)\right)}
{\lceil{n_j}^2\epsilon\rceil}d\mu (x)\right\vert\le \eta .
$$

So
$$
h_\mu^{2D} (F)=
\lim_{\epsilon\to 0}\liminf_{n_j\to\infty}\int_X \frac{-\log \left(\mu \left(\alpha_
{\lceil\epsilon\times n_j\rceil}^{n_j,F}(x)\right)
\right)}{\lceil\epsilon\times {n_j}^2\rceil}d\mu (x).
$$

As the sequences $\left(\int_X\frac{-\log\left( \alpha_p^n(x)\right)}{n}d\mu (x)\right)$  and
$\frac{1}{p_i}h_\mu (\alpha_{p_i},F)$ are a decreasing sequence then
$\int_X\frac{-\log\left( \alpha_{\lceil n\epsilon\rceil}^n(x)\right)}{\lceil n^2\epsilon\rceil}d\mu (x)
\le h_\mu^{2D}(F)$
then  we can conclude.

\hfill$\Box$

\begin{rem}
This second relation is not usefull to establish an inequality with the shift entropy and
the directional Lyapunov exponents.
\end{rem}

\subsection{First inequality}

\begin{pro}
For all  $F$ $\mu$ ergodic two-dimensional (CA),  one has
$$
h_\mu^{2D}(F)\le h_\mu (\sigma_1,\sigma_2)\times
(\lambda (0)+\lambda (\pi)+1)\times(\lambda (\pi /2)+\lambda (3\pi /2)+1).
$$
\end{pro}
{\it Proof}

Using Proposition \ref{firstinequality} we only need to show that

$$
\int_X \lim_{n\to\infty}\frac{-\log \mu (\alpha_n^{n,F}(x))}{n^2} d\mu (x)
$$
$$
\le (\lambda (0)+\lambda (\pi)+1)\times(\lambda (\pi /2)+\lambda (3\pi /2)+1).
$$

Let $R_n$ be the rectangular set of couple of integers $(i,j)$ such that

$-\Lambda_n^F(3\pi /2)\le i\le \Lambda^F_n(\pi /2)+n$ and
$\Lambda_n^F(\pi )\le i\le \Lambda^F_n(0)+n$.

From the definition of $T_n^*$, clearly one has $T_n^*\subset R_n$, so
using Proposition \ref{pro1} we obtain
$$
\alpha_n^{n,F}(x)\supset (\vee_{(i,j)\in T^*_n} \sigma_1^i\circ\sigma_2^j \alpha_n)(x)
\supset (\vee_{(i,j)\in R_n} \sigma_1^i\circ\sigma_2^j \alpha_n)(x).
$$

Put $(\vee_{(i,j)\in R_n} \sigma_1^i\circ\sigma_2^j \alpha_n)(x)=\alpha_n^R(x)$
 we have
$$
h_\mu^{2D}(F)\le\int_X \lim_{n\to\infty}\frac{-\log \mu\left(\alpha_n^R(x)\right)}{n^2} d\mu (x).
$$
So we can write
$$
h_\mu^{2D}(F)\le\int_X \lim_{n\to\infty}\frac{-\log \mu(\alpha_n^R(x))}{\vert R_n\vert}
\times\frac{\vert R_n\vert}{n^2} d\mu (x).
$$

From de definition of the two-sided shift one has
$$
\int_X\lim_{n\to\infty}\frac{-\log \mu(\alpha_n^R(x))}{\vert R_n\vert}d\mu (x)=
h_\mu (\sigma_1, \sigma_2,\alpha_1)=h_\mu (\sigma_1, \sigma_2).
$$

Hence we have
$$
h_\mu^{2D}(F)\le h_\mu (\sigma_1,\sigma_2)\times\lim_{n\to\infty}\frac{\vert R_n\vert}{n^2}.
$$

The surface $R_n$ is the surface of a rectangle of width $\Lambda^F_n(0)+\Lambda^F_n(\pi)$ and length
$\Lambda^F_n(\pi/2)+\Lambda^F_n(3\pi/2)$.
So $\vert R_n(x)\vert =(\Lambda^F_n(0)+\Lambda^F_n(\pi)+n)\times (\Lambda^F_n(\pi/2)+\Lambda^F_n(3\pi/2)+n)$.

Then
$$
\lim_{n\to\infty}\frac{\vert R_n\vert}{n^2}=
(\frac{\Lambda^F_n(0)}{n}+\frac{\Lambda^F_n(\pi)}{n}+\frac{n}{n})\times (\frac{\Lambda^F_n(\pi/2)}{n}+
\frac{\Lambda^F_n(3\pi/2)}{n}+\frac{n}{n})
$$

$$
=(\lambda (0)+\lambda (\pi)+1)\times(\lambda (\pi /2)+\lambda (3\pi /2)+1).
$$

So we can conclude.
\hfill$\Box$


\section{The continuity of the Lyapunov Exponents}


The study of the way the information propagate in a specific direction is interesting in itself.
The result in this section permit to establish a better upper bound for the value of the
always finite entropy because we now consider  the Lyapunov surface which is in general
leather than the simple squarre defined in the previous equality.

\begin{lem}
For all $n\in\NN^*$ big enought and for all  $\theta\in [0, 2\pi]$ there exits $\delta>0$
such that there exists a positive  integer $K$ such that for all $\theta'\le \delta$ one has
$$
\left\vert \frac{\Lambda^F_n(\theta +\theta')}{n}-\frac{\Lambda^F_n(\theta )}{n}\right\vert
\le K\times\vert\theta' \vert .
$$
\end{lem}
{\it Proof}

First recall that $\Lambda^F_n(\theta )=\max \{\Lambda_n^F(\theta )(x)\vert x\in X\}$; clearly we have
$$
\Lambda_n^F(\theta )=
\max\{s=\sqrt{i^2+j^2}\vert F^k(x)(a)=F^k(y) (a)\vert x\in X\; ;\; 0\le k\le n\}
$$
where $a=\left( (-r,r);(0,0)\right)$.

From de definition of $G_nF$, if  $(i,j)=G_nF(\theta )$, then for each $x\in X$ we obtain
$F^k(x)\left( (-r,r);(0,0) \right)=F^k(y)\left( (-r,r);(0,0) \right)$ where $y\in W^\theta_{(i,j)}(x)$.

Remark that there exist $\delta\in [0,2\pi]$ such that if $0\le \theta' \le\delta$ then
$\vert \tan{\theta'}\vert\le 2\theta')$.

As the information can not propagate more than $rn$ coordinates in $n$ iterations, it follows that
$F^k(x)\left( (-r,r);(0,0) \right)=F^k(y)\left( (-r,r);(0,0) \right)$ if
$y\in W^{\theta +\theta'}_{(u;v)}(x)$; $u=i+\lceil (r+1)n\times 2\vert\theta'\vert\times \cos{\theta}\rceil$ and
$v=j+\lceil (r+1)n\times 2\vert\theta'\vert\times \sin{\theta}\rceil$.

This implies that
$\Lambda_n^F(\theta +\theta')=\sqrt{u^2+v^2}=\sqrt{i^2+j^2}+(r+1)n\times 2\vert \theta'\vert$.


Finaly it follows that
$$
\left\vert \frac{\Lambda^F_n(\theta +\theta ')}{n}-\frac{\Lambda^F_n(\theta)}{n}\right\vert
\le 2(r+1)\vert\theta'\vert .
$$

\hfill$\Box$

\begin{pro}
For each positive integer $n$ the map $\frac{\Lambda^F_n(\theta)}{n}$ are continuous map,
moreover the map the sequence $\frac{\Lambda^F_n(\theta)}{n}$ converge uniformly to
 $\lambda (\theta )$ which is a continuous map.
\end{pro}
{\it Proof}

From Lemma 2 one has for all $n\in\NN^*$ and for all for all $\theta\in [0, 2\pi]$ there exits $\delta>0$
such that for all $\theta'\le \delta$ one has
$$
\left\vert \frac{\Lambda^F_n(\theta +\theta')}{n}-\frac{\Lambda^F_n(\theta )}{n}\right\vert \le K\times\theta' .
$$
So these maps are uniformly equicontinuous maps.

As the maps are equibounded by $r$ and are uniformly equicontinuous then we can use the
Ascoli-Arzela theorem which told us that that there exists a subsequence $n_i$ such that
the sequence $\frac{\Lambda^F_n(\theta )}{n}$ are uniformly convergent.
As from Lemma 1 the sequence $\frac{\Lambda^F_n(\theta )}{n}$ converge then we can conclude saying that
the maps $\frac{\Lambda^F_n(\theta )}{n}$ converge uniformely to $\lambda (\theta )$.

\hfill$\Box$


\section{The inequality}

\begin{thm}
For a shift ergodic measure $\mu$ of a two-dimensional shift and a two-dimensional cellular automata, we have

$h^{2D}_\mu (F) \le h_\mu (\sigma_1,\sigma_2)\times \left(\int_0^{2\pi}\frac{\lambda^2 (\theta )}{2}
+\sqrt {2}\int_0^{\pi}\lambda (\theta)) d\theta+\int_0^\pi d\theta\right)$.
\end{thm}
\paragraph{Proof}

From Proposition 2 one has
$$
h^{2D}_\mu (F) \le \int_X\liminf_{n\to\infty} \frac{-\log \mu \left(\alpha^{n,F}_{n}(x)\right)}{n^2}d\mu (x).
$$

From Proposition \ref{pro1} one has $\left(\vee_{(i,j)\in T_n^*}
\sigma_1^{i}\sigma_2^{j}(\alpha_1)\right)(x)=\alpha_n^{\sigma ,T_n^*}(x)\subset \alpha^{n,F}_n(x)$,
thus
$$
h^{2D}_\mu (F) \le \int_X\liminf_{n\to\infty} \frac{-\log \mu \left(\left(\vee_{(i,j)\in T_n^*}
\sigma_1^{i}\sigma_2^{j}(\alpha_1)\right)(x)\right)}{n^2} d\mu (x)
$$
and
$$
h^{2D}_\mu (F)\le \int_X \liminf_{n\to\infty}\frac{-\log\mu\left(\alpha_n^{\sigma ,T_n^*}(x)\right)}
{\vert T_n^*\vert}\times\frac{\vert T_n^*\vert}{n^2}d\mu (x).
$$

Since $(T^*_n)_{n\in\NN}$ is a special averaging sequence, we can use the extented Shannon-Breiman-McMillan
 Theorem (cf Orstein \cite{Or83}) which implies that for $\mu$-almost all points $x$ one has

$$
\lim_{n\to\infty}\frac{-\log\mu\left(\alpha_n^{\sigma ,T_n^*}(x)\right)}
{\vert T_n^*\vert}=h_\mu (\sigma_1,\sigma_2,\alpha_1)=h_\mu (\sigma_1,\sigma_2).
$$

Hence one obtains

$$
h^{2D}_\mu (F)\le h_\mu (\sigma_1,\sigma_2)\times  \liminf_{n\to\infty}
\frac{\vert T_n^*\vert}{n^2}.
$$

From Lemma 1 one has
$$
 \liminf_{n\to\infty}\frac{\vert T_n^*\vert}{n^2}=
\liminf_{n\to\infty}\frac{\vert T_n^{**}\vert}{T_n^{**}}\frac{T_n^*}{n^2}
$$

$$
=\liminf_{n\to\infty}\frac{\vert T_n^{**}\vert}{n^2}\le  \liminf_{n\to\infty}
\frac{\vert T_n\vert}{n^2}.
$$

Now we try to evaluate $\liminf_{n\to\infty}\frac{\vert T_n\vert}{n^2}$.
As $\vert T_n\vert$ is a sum of triangles of height $h_i=\Lambda^F_n(\frac{2\pi i}{n})$ and
base $b_i=\Lambda^F_n(\frac{2\pi i}{n})\times \tan (\frac{2\pi}{n})$ when
$\lfloor\frac{n}{2}\rfloor\le i\le n$.
For a real $\eta$ small enought and $n$ big enought we have
$$
\frac{\vert T_n\vert}{n^2}=\frac{\frac{2\pi}{n}\times\frac{1}{2}\left(
\sum_{i=0}^{\lceil\frac{n-1}{2}\rceil}\left(\Lambda^F_n(\frac{2\pi i}{n})+\sqrt{2}n\right)^2
+\sum_{i=\lfloor\frac{n}{2}\rfloor}^{n}\Lambda_n^2(\frac{2\pi i}{n}\right)+\eta}{n^2}
$$
\small
$$
=\frac{2\pi}{n}\left(\sum_{i=0}^{\lceil\frac{n-1}{2}\rceil}\left(\frac{1}{2}
\left(\frac{\Lambda^F_n(\frac{2\pi i}{n})}{n}\right)^2+
\sqrt(2)\frac{\Lambda^F_n(\frac{2\pi i}{n})}{n}+1\right)+\sum_{i=\lfloor\frac{n}{2}\rfloor}^{n}
\frac{1}{2}(\frac{\Lambda^F_n(\frac{2\pi i}{n})}{n})^2+\frac{\eta}{n}\right).
$$
\normalsize

From Proposition 5 the sequence of maps $\frac{\Lambda_n^F}{n}$ converge uniformely to the continuous map
$\lambda$. This implies that

$$
\lim_{n\to\infty}
\frac{2\pi}{n}\times \sum_{i=0}^{\lceil\frac{n-1}{2}\rceil}\left(\frac{1}{2}
\left(\frac{\Lambda^F_n(\frac{2\pi i}{n})}{n}\right)^2+\sqrt{2}\frac{\Lambda^F_n(\frac{2\pi i}{n})}{n}+1\right)
$$

$$
=\lim_{n\to\infty}
\frac{2\pi}{n}\times \sum_{i=0}^{\lceil\frac{n-1}{2}\rceil}\left(\frac{1}{2}
\left(\lambda (\frac{2\pi i}{n})\right)^2+\sqrt(2)\lambda (\frac{2\pi i}{n})+1\right).
$$

Therefore using the uniform continuity of the map $\lambda $ and the Riemman definition of the integral
we obtain $\liminf_{n\to\infty}\frac{\vert T_n(x)\vert}{n^2}=\lim_{n\to\infty}\frac{\vert T_n(x)\vert}{n^2}$ and

$$
\lim_{n\to\infty}\frac{\vert T_n(x)\vert}{n^2}=\int_\pi^{2\pi}\frac{\lambda^2 (\theta )}{2}d\theta +
\int_0^{\pi} \frac{\lambda^2 (\theta )}{2}d\theta
+\sqrt (2)\int_0^{\pi}\lambda (\theta ) d\theta+\int_0^{\pi}d\theta .
$$

Finaly
\small
$$
h^{2D}_\mu (F)\le
h_\mu (\sigma_1,\sigma_2)\times
\left(\int_0^{2\pi} \frac{\lambda^2 (\theta )}{2}d\theta
+\sqrt (2)\int_0^{\pi}\lambda (\theta) d\theta+\int_0^{\pi}d\theta\right).
$$
\normalsize

\hfill$\Box$

\begin{rem}
We could obtain a better upper bound if we will be able to estimate the value of
$\liminf_{n\to\infty}\frac{\vert T^*_n\vert}{n^2}$ instead of
$\liminf_{n\to\infty}\frac{\vert T_n\vert}{n^2}$.
\end{rem}


\subsection{Computation of $\lambda (\theta )$}

We compare here the two kind of upper bound (Proposition \ref{firstinequality} and Theorem 1)
in the case of the example $F_3$ and different extensions of this map.

In order to compute the first upper bound, remark that $\lambda(0)=1$, $\lambda(\pi )=0$;
$\lambda(\pi/2 )=1$
and $\lambda (3\pi/2) =0$.

So
$\lambda_R^{F_3}=\left(\lambda (0)+\lambda (\pi)+1\right)$ $\times\left(\lambda (\pi /2)+\lambda (3\pi /2)+1\right)=4.$
 We obtain the inequality
$$
h_\mu^{2D}(F)=2\log (2)\le 4h_\mu(\sigma_1,\sigma_2)=4\log (2).
$$

Now try to evaluate the second upper bound :
$$
\lambda_T^{F_3}=\int_0^{2\pi}\frac{(\lambda (\theta ))^2}{2}d\theta +\sqrt (2)\int_0^{\pi} \lambda (\theta )d\theta
 +\int_0^{\pi}d\theta .
$$

As no perturbation came from the left or under part of the latttice we only consider the Lyapunov exponents
between $0$ and $\pi/2$, that it to say we use $S_n=R_n\cap T_n\subset T_n^{**}$ where $R_n$ is the set defined
in the proof of Proposition \ref{firstinequality}.
So we have to compute
$$
\int_0^{\pi/2}\frac{(\lambda (\theta ))^2}{2}d\theta +\sqrt (2)\int_0^{\pi/2} \lambda (\theta )d\theta
 +\frac{\pi}{2}.
$$

If $\theta\in [0,\pi/4]$ then $\lambda (\theta )=\cos (\theta )$ and if
$\theta\in [\pi/4,\pi/2]$ then $\lambda (\theta )=\cos (\pi/2-\theta )$.

Then $\lambda_T^{F_3}=\int_0^{\pi/2}\frac{1+\cos(\theta )}{2}d\theta +\sqrt{2}\int_0^{\pi/4}
\cos(\theta)d\theta+\pi/2=
(\frac{\pi}{8}+\frac{1}{4})+\sqrt{2}+\frac{\pi}{2}$.
We have $\lambda_T^{F_3}=3.628\le 4$. More generaly for any integer $k\ge 1$ denote by
$F'_k$ a (CA) defined by the local rule $[F'_k(x)]_{(i,j)}=x_{(i+k,j)}+x_{(i,j+k)}$ mod 2.
In this case we have $\lambda_R^{F_3}=k^2+2k+1$ and
$\lambda_T^{F_3}=k^2(\frac{\pi}{8}+\frac{1}{4})+\sqrt{2}k+\frac{\pi}{2}$.
The difference between the two exponents increase when the value of $k$.

\subsection{Conclusion}
Because the value of the standard entropy appears to be zero or not finate for a two-dimensional (CA),
 the Always Finite Entropy is a map which allows to extend two kinds of results which appear in the
one-dimensional case. The first one is the value of metric entropy of a (CA) for additive rules with the
uniform measure. The second one is to connect the value of the entropy of the (CA) with the entropy of
the shift and the speed of propagation of information (Lyapunov exponents).

Is this version of metric entropy will be useful in more general dynamical systems for which the
value of the entropy is not finate.

\subsection{Questions}
In these following questions $F$ always represent a two-dimensional cellular automaton.

\begin{itemize}

\item If $F$ is bijective is it true that $h(F)<\infty$ and $h^{2D}(F)=0$?

\item Is it possible that $h_\mu (F)>0$ but $h_\mu (F)<\infty$ (Shereshevsky conjecture)?

\item Is it possible that $h_\mu^{2D}(F)=0$ and $h_\mu (F)>0$ ?

\item Is it possible to define directional exponents which do not depend on the maximum
and nevertheless are continuous with the direction $\theta$ ?

\item Is it possible to find a better upper bound which is equal to zero when
the cellular automaton is the identity map?

\item In the three examples the value of $h_\mu^{2D}(F)$ depends linearly of the radius $r$ of the local
rule. The expression
of the upper bound depend on some square of  $\lambda(\theta )$. We can wonder if there exist some
cellular automaton with the property that  $h_\mu^{2D}$ is proportional of the square of the radius
of the local rule.




\end{itemize}


\end{document}